%% file: six_pile_remoteness.tex
\documentclass[11pt]{amsart}
\usepackage[T1]{fontenc}
\usepackage[utf8]{inputenc}
\usepackage{lmodern}
\usepackage[margin=1in]{geometry}
\usepackage{amsmath,amssymb,amsthm}
\usepackage{booktabs,array}
\usepackage{needspace}
\usepackage{microtype}
\usepackage{xurl}
\usepackage[hidelinks]{hyperref}
\hypersetup{
  pdftitle={A general reduction and explicit remoteness formulas for exact slow 2-Nim},
  pdfauthor={Ra\'ul S\'anchez Gal\'an and Yijie Xu},
  pdfkeywords={exact slow Nim, Smith remoteness, subtraction games, GM rule}
}

\newtheorem{theorem}{Theorem}[section]
\newtheorem{lemma}[theorem]{Lemma}
\newtheorem{proposition}[theorem]{Proposition}
\newtheorem{corollary}[theorem]{Corollary}
\theoremstyle{remark}

\newtheorem{example}[theorem]{Example}
\numberwithin{equation}{section}
\newcommand{\Z}{\mathbb Z}
\newcommand{\one}{\mathbf 1}
\newcommand{\parround}{\Pi}
\newcommand{\ind}[1]{\left[#1\right]}
\title[Remoteness in exact slow 2-Nim]{A general reduction and explicit remoteness formulas for exact slow $2$-Nim}
\author[Ra\'ul S\'anchez Gal\'an]{Ra\'ul S\'anchez Gal\'an\textsuperscript{1}}
\author[Yijie Xu]{Yijie Xu\textsuperscript{2}}
\thanks{\textsuperscript{1}\,4i Intelligent Insights, 41092 Sevilla, Spain.
  \href{mailto:r.sanchez@4i.ai}{\nolinkurl{r.sanchez@4i.ai}}}
\thanks{\textsuperscript{2}\,The Hong Kong University of Science and Technology,
  Guangzhou, China.
  \href{mailto:yxu409@connect.hkust-gz.edu.cn}{\nolinkurl{yxu409@connect.hkust-gz.edu.cn}}}
\date{}

\begin{document}
\begin{abstract}
In exact slow $2$-Nim, a move removes one token from each of two nonempty
piles. For any number of piles, we relate Smith remoteness in the exact
game to the remoteness $D$ of Moore's slow game on the piles that remain
after a largest pile is removed. If the total number of tokens $N$ is odd,
the exact game has remoteness $\min(D,\lfloor N/2\rfloor)$. If $N$ is even,
its remoteness equals $D$ when $2D\le N$ and is either $N/2$ or $N/2-1$
otherwise. Both values can occur for the same $N$, $D$, and number of
piles. We determine $D$ explicitly for five piles. This gives a six-pile
remoteness formula and a rule for selecting an optimal move. Both take
$O(b)$ bit operations when each pile size has at most $b$ bits. We also
give a closed formula for positions with any number of unit piles and
two further piles of arbitrary size. For the GM move rule, we prove the
proposed descriptions of the exceptions in the four- and five-pile games
and give the exact threshold at which a specified six-pile move ceases
to be optimal as the largest pile grows. Finally, we classify all
outcome-preserving sorted translations for five and six piles.
\end{abstract}
\maketitle

\section{Introduction}\label{sec:introduction}

Knowing which player wins an impartial game does not determine how
long optimal play lasts. The parity of Smith's remoteness function
distinguishes winning and losing positions. Its value gives the
duration of play when the winner seeks the earliest possible
victory and the loser delays defeat as long as possible
\cite{Smith1966}. Explicit remoteness formulas therefore refine an
outcome classification and distinguish winning moves by how quickly
they force a win. For Moore's Nim, Wythoff's game and Euclid's
game, Boros, Gurvich, Makino and Vyalyi studied the complexity of
computing remoteness \cite{Boros2024}.

We study the normal-play game of \emph{exact slow $2$-Nim}. A position
consists of piles of tokens, and each move removes exactly one token
from each of exactly two distinct nonempty piles. A player with no
legal move loses. We write $R_n$ for the remoteness function with $n$
piles. Empty piles may be inserted or deleted without changing the
game.

Exact slow Nim has been studied through its outcome classes,
Sprague--Grundy function and move-selection rules. Chikin, Gurvich,
Knop, Paterson and Vyalyi \cite{Chikin2021} characterized the
P-positions for five and six piles and proved several
Sprague--Grundy shift identities. Gurvich, Martynov, Maximchuk and
Vyalyi \cite{GurvichEtAl2024} determined remoteness in exact slow
$k$-Nim with $k+1$ piles; their results include the three-pile game
considered here. More recently, Gurvich, Maximchuk and Naumova
\cite[Sections~2.2--2.3]{GurvichMaximchukNaumova2026} proposed explicit
patterns for the four- and five-pile positions at which the GM rule
fails to reduce remoteness by one. Their formulas describe the exceptions they found computationally.

Gurvich, Heubach, Ho and Chikin \cite[Section~4]{Slow2020}
used the correspondence between exact two-pile moves and moves in
Moore's slow game on one fewer pile. Theorem~\ref{thm:transfer}
turns this correspondence into a formula for remoteness at every
odd total and a sharp two-value alternative at even totals.
It also identifies the least value of the largest pile for which the
two remotenesses agree. The proof uses a bound on the sum of the
Moore piles plus their maximum, which controls the moves lost when
the two largest exact-game piles are tied.

For six piles, the known outcome classification resolves the two-value
alternative. It then suffices to evaluate remoteness in Moore's slow
$2$-Nim on five piles. Its P-positions were determined in
\cite{Slow2020}. Three parity bits
determine how each move changes two linear bounds.
The calculation gives both the six-pile value
(Theorem~\ref{thm:six}) and an optimal move
(Proposition~\ref{prop:direct-six}).

For five piles, one linear bound suffices. Comparison with
the half-total gives a short description of every exception to
the deterministic GM rule and an optimal replacement move.
Corollary~\ref{cor:all-gm} also determines when every admissible
choice of odd pivots fails. For six piles, the behavior changes.
Even at a losing position with all six piles even, the GM move can
bring defeat two moves closer. Corollary~\ref{cor:even-six} describes
exactly when this happens. More generally, Theorem~\ref{thm:pair-threshold}
gives a threshold for each pair as the largest pile varies.

Chikin et al.~\cite{Chikin2021} exhibited translations preserving
the P-positions. Theorem~\ref{thm:all-shifts} classifies all such
translations of sorted positions with five or six piles.
Every one increases remoteness by half the total number of added
tokens. For general vector subtraction games, a separate two-step
criterion shows that remoteness increases by two whenever the subtraction set
satisfies $S=c-S$. Applied to a position represented by the number
of unit piles and the sizes of two further piles, it gives
an explicit remoteness formula with no bound on the pile count
(Corollary~\ref{cor:two-large}).

\section{Remoteness and sorted moves}\label{sec:preliminaries}

All games below are impartial and use the normal-play convention.
Every position has finitely many options, and every sequence of moves
from a given position is finite. A
\emph{P-position} is losing for the player to move, and an
\emph{N-position} is winning for that player. Smith remoteness is
defined recursively by
\begin{equation}\label{eq:recurrence}
 R(x)=
 \begin{cases}
 0,&x\text{ is terminal},\\[2pt]
 1+\min\{R(y):x\to y,\ R(y)\text{ even}\},
   &\text{if an even-remoteness option exists},\\[2pt]
 1+\max\{R(y):x\to y\},&\text{otherwise}.
 \end{cases}
\end{equation}
It follows by induction that $R(x)$ is even precisely at P-positions.
A move $x\to y$ is \emph{remoteness-optimal} if
$R(y)=R(x)-1$. The following form of the recurrence will be used
throughout.

\begin{lemma}[Verification criterion]\label{lem:verification}
Let $F$ be a nonnegative integer-valued function on the positions of
a game. Suppose that $F$ vanishes at terminal positions and, at every
nonterminal position $x$, the following conditions hold.
\begin{enumerate}
\item If $F(x)$ is even, every option has an odd $F$-value at most
      $F(x)-1$, and some option has value $F(x)-1$.
\item If $F(x)$ is odd, some option has value $F(x)-1$, and every
      option with even $F$-value has value at least $F(x)-1$.
\end{enumerate}
Then $F=R$.
\end{lemma}

\begin{proof}
Induct on the maximum length of play from $x$. In the first case,
the maximum option value is $F(x)-1$ and all options have odd
remoteness. In the second case, the minimum even option value is
$F(x)-1$. These are exactly the two nonterminal clauses of
\eqref{eq:recurrence}.
\end{proof}

We always write pile sizes in nondecreasing order. Since the piles
are interchangeable, this loses no information about remoteness.
When equal piles occur, it is useful to choose moves that preserve
the order without a subsequent permutation.

\begin{lemma}[Canonical moves]\label{lem:canonical}
Suppose $x_1\le\cdots\le x_n$ and a move decreases a prescribed
number of positive coordinates by one. Every sorted successor has a
representative $x-\epsilon$, where $\epsilon_i\in\{0,1\}$ and
\begin{equation}\label{eq:canonical}
 x_i=x_{i+1}\quad\Longrightarrow\quad
 \epsilon_i\ge\epsilon_{i+1}.
\end{equation}
For this representative, $x-\epsilon$ is already sorted.
\end{lemma}

\begin{proof}
Within each block of equal coordinates, move all selected indices to
the beginning of the block. This does not change the multiset of
successor pile sizes. It gives \eqref{eq:canonical}. Between two
different blocks, the original sizes differ by at least one, so a
unit decrement cannot reverse their order.
\end{proof}

We call such a representative a \emph{canonical move}. Primes will
denote quantities evaluated at its successor. For a statement $P$,
the indicator $[P]$ is one if $P$ holds and zero otherwise. Unless a
different modulus is displayed, congruences in parity conditions
are modulo two.

For integers $N\ge0$ and $D$, define the parity adjustment
\begin{equation}\label{eq:pi}
 \parround(N,D)=\left\lfloor\frac N2\right\rfloor
 -\ind{N\text{ even and }\lfloor N/2\rfloor\not\equiv D\pmod2}.
\end{equation}
If $N$ is even, this is the largest integer at most $N/2$ with the
parity of $D$. If $N$ is odd, it is simply $\lfloor N/2\rfloor$.
Two identities, immediate from the definition, are
\begin{align}
 \parround(N-2,D-1)&=\parround(N,D)-1 \quad (N\ge2),
     \label{eq:pi-minus}\\
 \parround(N+4,D+2)&=\parround(N,D)+2.
     \label{eq:pi-plus}
\end{align}

\section{Reduction to Moore's slow game}\label{sec:transfer}

Let $x=(z,f)$ be a sorted position with $n\ge2$ piles, where $f$
is its largest entry and $z$ consists of the other $n-1$ piles.
Write $S=\sum z_i$, $m=\max z_i$, $N=S+f$ and $H=\lfloor N/2\rfloor$.
In Moore's slow $2$-Nim, a move selects one or two nonempty piles
and removes one token from each. Write $D(z)$ for its remoteness.
An exact move always projects to a Moore move after a largest pile
is deleted. This correspondence was used for outcome and
Sprague--Grundy calculations in \cite[Section~4]{Slow2020}.
The following theorem determines how much remoteness it preserves.

\begin{theorem}[Largest-pile reduction]\label{thm:transfer}
Put $D=D(z)$. For any number of piles,
\begin{align}
 N\text{ odd}&\quad\Longrightarrow\quad R_n(z,f)=\min(D,H),\label{eq:odd-transfer}\\
 N\text{ even},\ D\le H&\quad\Longrightarrow\quad R_n(z,f)=D,\label{eq:even-transfer}\\
 N\text{ even},\ D> H&\quad\Longrightarrow\quad R_n(z,f)\in\{H-1,H\}.\label{eq:even-band}
\end{align}
The least $f\ge m$ for which $R_n(z,f)=D$ is therefore
\begin{equation}\label{eq:saturation}
               f_* = \max\{m,\ 2D-S\},
   \qquad R_n(z,f)=D\quad\Longleftrightarrow\quad f\ge f_* .
\end{equation}
\end{theorem}

\begin{lemma}[A lower bound for Moore remoteness]\label{lem:moore-bound}
Let $T(z)=\sum z_i+\max z_i$. Then
\[
                  T(z)\le 2D(z)+[D(z)\text{ odd}].
\]
In particular, $D(z)\ge\lfloor T(z)/2\rfloor$.
\end{lemma}
\begin{proof}
If at most one pile is positive, then $D=m$ and $T=2m$.
Otherwise, induct on $S$. If $D$ is odd, choose an even option
of value $D-1$. A move decreases the sum by at most two and the
maximum by at most one. Thus $T\le 2(D-1)+3$.
If $D$ is even, leave a designated largest pile untouched and
reduce one other positive pile. This reduces $T$ by one. The
successor has odd remoteness at most $D-1$, giving
$T\le 2(D-1)+1+1$.
\end{proof}

\begin{proof}[Proof of Theorem~\ref{thm:transfer}]
By Lemma~\ref{lem:canonical}, an exact move has a sorted
representative. Deleting its last coordinate leaves one or two
decrements in $z$. Every Moore move on two piles lifts:
decrease those two and leave $f$ untouched. When $f>m$, every
single decrement also lifts by using $f$ as its partner.
When $f=m$, a single decrement of a unique largest pile of $z$
still lifts, since both copies of that size can be decreased.
Other single decrements can fail to lift at a tie. For instance, all
exact moves on $(1,1,1)$ project to a move on both remaining piles.
The exact position is terminal if and only if its projection is terminal.

For a list of option values, write $B$ for its preferred member in
Smith's recurrence. An even value is preferred to any odd value.
Among even values, smaller values are preferred. Among odd values,
larger values are preferred. Truncation at $t$, given by
$C_t(v)=\min(v,t)$, preserves this weak order. Indeed, the values
at least $t$ lie between the even values below $t$ and the odd
values below $t$ in that order. Hence truncation commutes with $B$,
and $C_t(1+v)=1+C_{t-1}(v)$.

Set $h=H$ for odd $N$ and $h=H-1$ for even $N$. We claim
\begin{equation}\label{eq:capped-transfer}
                     C_h(R_n(x))=C_h(D(z)).
\end{equation}
At a terminal position, both values before truncation are zero.
For $N\le2$, the only nonterminal position has two unit piles
and both values are one. Induct on $N$ at all other positions.
For each exact option, the truncation level is $h-1$. If $f>m$, the projected
options exhaust the Moore options. Apply the induction hypothesis
to their truncated values, take their preferred member and add one.
The two identities for $C_h$ give \eqref{eq:capped-transfer}.

If $f=m$, then $T(z)=N$ and Lemma~\ref{lem:moore-bound} gives
$D(z)\ge H$. After any exact move, the sum of the projected pile
sizes decreases by at most two and their maximum by at most one. Thus
$T(z')\ge N-3$ and
$D(z')\ge\lfloor(N-3)/2\rfloor=h-1$.
The induction hypothesis gives $R_n(x')\ge h-1$ for every option.
Consequently $R_n(x)\ge h$, proving \eqref{eq:capped-transfer}
also at a tie.

Every exact play has length at most $H$. For odd $N$ this proves
\eqref{eq:odd-transfer}. For even $N$, \eqref{eq:capped-transfer}
already gives \eqref{eq:even-band} and exact equality when $D<H-1$.
It remains to show that the two games have the same winner when
$N$ is even and $D\le H$.

Follow the winning player's optimal strategy in the projected game.
On each of the winner's turns, the Moore value drops from odd $D$ to even
$D-1$. On the other player's turns, every projected move gives an odd
value at most $D-1$. The inequality $D\le H$ is therefore preserved
as long as the winner's moves can be lifted. Only a tied top needs
attention. There $D\ge H$ by the lemma, so $D=H$. If an optimal
single decrement left the projected maximum unchanged, its successor
would have $T(z')=2H-1$. Its even value $D-1$ would give
$T(z')\le 2H-2$, a contradiction. Thus an optimal single decrement
lowers a unique projected maximum and lifts. All two-pile moves
lift as well. The winner can carry out the strategy to a common
terminal position, so $R_n(x)$ and $D$ have the same parity.
Equation~\eqref{eq:capped-transfer} now gives $R_n(x)=D$ for
$D=H-1$ and $D=H$. Finally, equality with $D$ holds exactly
when $N\ge 2D$, which is \eqref{eq:saturation}.
\end{proof}

\begin{corollary}[Optimal moves under projection]\label{cor:tied}
Suppose $x$ is nonterminal. If $f>m$ and either $N$ is odd or
$D\le H$, every optimal Moore move lifts to an optimal exact move.
If $f=m$ and $N$ is odd, every move
is remoteness-optimal. If $f=m$ and $N$ is even, every move is
optimal at a P-position, and every move to a P-position is optimal
at an N-position.
\end{corollary}
\begin{proof}
When the largest pile is unique, an optimal Moore successor has value
$D-1$ and the half-total drops to $H-1$. The first two lines of
Theorem~\ref{thm:transfer} give the claimed decrement.
At a tie the preceding proof gives $h-1\le R_n(x')\le H-1$ at every option.
For odd $N$ these bounds coincide. For even $N$ they are consecutive
integers, so all options of the same outcome have the same remoteness.
Apply \eqref{eq:recurrence}.
\end{proof}

The ambiguity in \eqref{eq:even-band} is sharp even when $n,N,D$
are all specified. Consider eight unit piles and $(2,3,3)$ padded
with five zeros. Both have $N=8$ and $D=5$. In the first example
the projected game removes one or two units. Its value on $3q+j$
units is $2q+[j>0]$ for $0\le j\le2$, by induction on the number
of units. Seven units thus give five, and the eight-pile exact game
lasts four moves. For two Moore piles $(a,b)$ the same recurrence
gives $a+b-[a,b\text{ both odd}]$. An even successor requires
decreasing every odd pile, and from an all-even position a single
decrement attains the largest odd value. This gives $D(2,3)=5$.
In the exact game, decreasing both threes
leaves $(2,2,2)$, whose only sorted option $(1,1,2)$ has remoteness
one. Thus $(2,3,3)$ can be won in three moves, and no move wins
immediately.

\section{Six piles}\label{sec:six}

Write $x=(a,b,c,d,e,f)$ in nondecreasing order and let $N$ be its
total. Define
\begin{align}
 A&=b+d+e, & B&=a+b+2d, & h&=\ind{c\not\equiv d},\label{eq:six-forms}\\
 E(x)&:\quad a\equiv b,\quad c\equiv d,\quad e\equiv a+d.&&&
 \notag
\end{align}
Put $p=0$ when $E(x)$ holds and $p=1$ otherwise. Set
\begin{equation}\label{eq:six-round}
 M=\max\{A,B-h\},\qquad D_6=M-\ind{M\not\equiv p}.
\end{equation}
Thus $D_6$ is the largest integer at most $M$ with parity $p$.

\begin{theorem}[Six-pile remoteness]\label{thm:six}
For every $0\le a\le b\le c\le d\le e\le f$,
\begin{equation}\label{eq:six}
 \boxed{\displaystyle R_6(a,b,c,d,e,f)=
   \min\{D_6,\parround(N,D_6)\}.}
\end{equation}
\end{theorem}

The two linear forms in \eqref{eq:six-forms} can be read from the
outcome theorem of Chikin et al.~\cite[Theorems~1--2]{Chikin2021}.
Their two boundary expressions are $N-2A$ and $N-2B+1$.
When $B>A$ and $N=2B-1$, the first is positive and the second vanishes.
Their criterion gives a P-position there
exactly when $N\equiv1\pmod4$ and $h=0$. In particular, if $B$
is odd and $h=1$, the half-total bound $B-1$ is even but the
position is winning. The correction in \eqref{eq:six-round} lowers
the odd cap to $B-2$. For example, at the boundary position $(1,2,2,3,3,6)$ we
have $A=8$, $B=9$, $h=1$ and $p=1$. Formula~\eqref{eq:six} gives
seven moves. Omitting $h$ would give eight and the wrong winner.

In particular, put $U=\min\{\lfloor N/2\rfloor,M\}$. Since each
term of the minimum in \eqref{eq:six} is lowered by at most one
from its corresponding bound,
\[
                     R_6(x)\in\{U-1,U\}.
\]
The duration is the even member of this pair at a P-position and the
odd member at an N-position. If each pile size has at most $b_0$ bits,
both the formula and the move rule below use $O(b_0)$ bit operations.

\subsection{The parity calculation and an optimal move}
\label{subsec:constructive}

Three bits record how a move changes the two bounds. Write
\[
 r=B-A=a+d-e,\qquad
 u=(b-a,\ d-c,\ r)\bmod 2=(i,j,k).
\]
We use strings such as $101$ for these triples and $\oplus$ for
componentwise addition modulo two. Then $p=[u\ne000]$ and $h=j$.
For a bit $v$, let $\rho_v(z)=z-((z-v)\bmod2)$, the largest integer
at most $z$ with parity $v$. Rounding commutes with the maximum. Since
$B\equiv i$ and $A\equiv i-k$, we can write
\begin{equation}\label{eq:offset-form}
 D_6=\max(X,Y),\qquad X=A-\alpha(u),\quad Y=B-\beta(u),
\end{equation}
where $\alpha(u)=(i-k-p)\bmod2$ and
$\beta(u)=j+((i-j-p)\bmod2)$. In particular
$\alpha(000)=\beta(000)=0$.

Number the piles from 1 to 6. Decreasing one coordinate changes $u$
by the corresponding column of
\[
                  (101,\ 100,\ 010,\ 011,\ 001,\ 000).
\]
A move on a pair adds the two columns. Table~\ref{tab:six-pairs}
groups the fifteen pairs by this sum $v$. The integers following
each pair are its losses $m=A-A'$ and $n=B-B'$. They are obtained
by adding the corresponding entries of $(0,1,0,1,1,0)$ for $A$
and $(1,1,0,2,0,0)$ for $B$.
For instance, pair $46$ changes $u$ by $011$ and lowers $A,B$ by
one and two, respectively.

\begin{table}[htbp]
\centering
\small
\renewcommand{\arraystretch}{1.18}
\begin{tabular}{@{}cccll@{}}
\toprule
$v$ & $\alpha(v)$ & $\beta(v)$ & Pairs and their $(m,n)$
   & Choice from $u=v$\\
\midrule
$001$ & 0 & 1 & $12,34\ (1,2)\quad 56\ (1,0)$
 & $12$ if $a>0$, else $56$\\
$010$ & 1 & 1 & $36\ (0,0)\quad 45\ (2,2)$ & $45$\\
$011$ & 0 & 1 & $35\ (1,0)\quad 46\ (1,2)$ & $46$\\
$100$ & 0 & 0 & $15,26\ (1,1)$ & $26$\\
$101$ & 1 & 0 & $16\ (0,1)\quad 25\ (2,1)$
 & $25$ if $e>d$, else $16$\\
$110$ & 0 & 2 & $14\ (1,3)\quad 23\ (1,1)$
 & $14$ if $a>0$, else $23$\\
$111$ & 1 & 2 & $13\ (0,1)\quad 24\ (2,3)$ & $24$\\
\bottomrule
\end{tabular}
\caption{Parity changes and losses of the two linear forms.
The last column is used in part~\textup{(iii)} of
Proposition~\ref{prop:direct-six}.}\label{tab:six-pairs}
\end{table}

\begin{proposition}[An optimal first move]\label{prop:direct-six}
At a nonterminal position, follow the first instruction that applies.
\begin{enumerate}
\item[(i)] If $f=e$ and either $N$ is odd or $u=000$, use any
legal canonical pair, such as the first two positive piles.
\item[(ii)] If $u=000$, decrease the first positive pile among
$1,\ldots,5$ and pile $6$.
\item[(iii)] Otherwise, use the last column of Table~\ref{tab:six-pairs}.
\end{enumerate}
The selected pair is canonical and reduces remoteness by one.
\end{proposition}

\begin{proof}[Proof of Theorem~\ref{thm:six} and Proposition~\ref{prop:direct-six}]
First consider Moore's slow game on $(a,b,c,d,e)$. Interpret a
pair containing index~6 as a single decrement of its other pile.
We verify that $D_6=\max(X,Y)$ is its remoteness.

The parity columns are distinct. From $u=000$, every option has
a nonzero triple $v$ and an odd value of $D'_6$. The losses in $(X,Y)$ are
$(m+\alpha(v),n+\beta(v))$, both at least one throughout
Table~\ref{tab:six-pairs}. Thus $D'_6\le D_6-1$. From $u\ne000$ an
even $D'_6$ requires $u'=000$, hence $v=u$. The losses are now
$(m-\alpha(v),n-\beta(v))$, both at most one, giving $D'_6\ge D_6-1$.

To attain the bound from $000$, decrease the first positive pile. Each
pair with index~6 lowers $X$ by one. If $a>0$, the chosen pair
$16$ also lowers $Y$ by one. If $a=0$, then $r=d-e\le0$ and
$D_6=X$. The upper bound on $D'_6$ therefore gives equality.
From a nonzero triple, the move in the last column returns to $000$. Both
losses are one except in three cases. In row $001$ with $a=0$,
$r=d-e$ is negative and odd, so $X\ge Y+2$ and pair $56$ has
losses $(1,-1)$. In row $101$ with $e=d$, $r=a$ is positive
and odd, so $Y\ge X+2$ and pair $16$ has losses $(-1,1)$.
In row $110$ with $a=0$, $r=d-e\le0$, so $X\ge Y+2$ and
pair $23$ has losses $(1,-1)$. The maximum drops by one in each case.
All choices respect the order. In particular, $i=1$ gives $b>a$,
$j=1$ gives $d>c$, and the pairs selecting pile~5 either
also select pile~4 or have $e>d$. The explicit tests in the
table ensure that the selected piles are positive.
When $e=0$ the value is zero, and otherwise $M\ge A\ge1$.
Lemma~\ref{lem:verification} proves the Moore formula.

For the exact game with even total, the P-positions are precisely
those with $u=000$ \cite[Theorem~1]{Chikin2021}. One can also read
this from the table. No move preserves $000$, and the last column
returns every nonzero triple to $000$. Indeed, a pair using pile~6
without pile~5 is chosen only in a row with $i\ne j$. The identity
$N\equiv(f-e)+i+j$ then gives $f>e$ for even $N$.
Terminal positions have triple $000$, so induction proves the
claim. Thus the exact outcome has the parity of $D_6$ at even
totals. Theorem~\ref{thm:transfer} now gives
\eqref{eq:six} for both total parities.

When $f>e$, the Moore move above lifts and has $D'_6=D_6-1$.
Equations~\eqref{eq:six} and \eqref{eq:pi-minus} show that it
is optimal in the exact game. At a tied top,
Corollary~\ref{cor:tied} allows any move if $N$ is odd or $u=000$.
In the remaining case the last column reaches a P-position, so
the same corollary proves optimality.
\end{proof}

\noindent\emph{Moves that can fail.}
Suppose each attempted move succeeds independently with probability
$q\in(0,1)$, and a failure passes the turn without changing the piles.
Put $\lambda=q/(2-q)$. By S\'anchez Gal\'an and Xu
\cite[Theorem~5.2]{SanchezGalanXu2026}, the player to move at $x$ wins
with probability $\tfrac12(1-(-\lambda)^{R_6(x)})$ under optimal play,
and an attempt is optimal exactly when the corresponding move lowers
$R_6$ by one. Thus Theorem~\ref{thm:six} gives the winning probability
explicitly, and the pair selected by Proposition~\ref{prop:direct-six}
is optimal for every such $q$. Theorem~\ref{thm:pair-threshold} determines
when attempting a specified pair maximizes this probability.

\begin{corollary}[Five-pile Moore remoteness]\label{cor:moore}
In Moore's slow $2$-Nim, the remoteness of $(a,b,c,d,e)$ is the
value $D_6$ defined in \eqref{eq:six-round}. At a nonterminal position, apply
parts~\textup{(ii)} and \textup{(iii)} of Proposition~\ref{prop:direct-six}
and omit index~6 from the selected pair. This gives an optimal move.
\end{corollary}
\begin{proof}
The preceding proof gives the formula for the projected game and
shows that the specified moves attain it.
\end{proof}

The even-total parity criterion changes for seven piles. For example,
the positions $x=(1,2,2,2,3,3,3)$ and $y=(1,2,2,4,5,5,5)$ have the
same ordered parities and both totals are divisible by four. The first
is losing and the second is winning. To see
this, after any move from $x$ other than a move on two threes,
decrease two untouched piles so that the four original sizes chosen
over the two moves are $1,2,2,3$. The result, with its zero omitted,
is $(1,1,2,2,3,3)$. After a move on two threes, decrease the one and
the remaining three, leaving six twos. Both replies reach a six-pile
P-position. From $y$, decreasing the one and the four leaves the
six-pile P-position $(2,2,3,5,5,5)$. Corollary~\ref{cor:tied}
and the bounds at a tied top give $R_7(x)=8$ and $R_7(y)=11$.

\section{The five-pile formula}\label{sec:formula}

For a sorted position $x=(a,b,c,d,e)$, put
\begin{equation}\label{eq:parameters}
 N=a+b+c+d+e,\qquad L=a+c+d,\qquad
 \Delta=N-2L=b+e-a-c-d.
\end{equation}
The parity condition
\begin{equation}\label{eq:q}
 Q(x):\qquad a\text{ even and }b\equiv c\equiv d
\end{equation}
determines the adjustment of $L$. Set
\begin{equation}\label{eq:chi}
 \chi=[L\text{ even and }Q(x)\text{ fails}],
 \qquad D=L-\chi.
\end{equation}
Equivalently, $D$ is the largest integer at most $L$ that is even
when $Q(x)$ holds and odd otherwise. Thus $D$ is a parity rounding
of the single linear quantity $L$.

\begin{corollary}[Five-pile remoteness]\label{thm:five}
For every $0\le a\le b\le c\le d\le e$,
\begin{equation}\label{eq:five}
 \boxed{\displaystyle
 R_5(a,b,c,d,e)=\min\{D,\parround(N,D)\}.}
\end{equation}
\end{corollary}
\begin{proof}
Insert a leading zero. The six-pile forms become $A=L$ and
$B=a+2c\le L$, and $E(0,a,b,c,d,e)$ is exactly $Q(a,b,c,d,e)$.
The rounding in \eqref{eq:six-round} gives $D=L-\chi$, so the
result follows from Theorem~\ref{thm:six}.
\end{proof}

For example, at $x=(2,3,4,4,8)$, we have $N=21$ and $L=10$,
and $Q(x)$ does not hold. Hence $D=9$ and $\parround(N,D)=10$, giving
$R_5(2,3,4,4,8)=9$.

The pair $(b,e)$ and the triple $(a,c,d)$ have totals differing
by $\Delta$. Its sign determines which term of \eqref{eq:five}
controls the value. Write $F(x)=R_5(x)$.

\begin{lemma}\label{lem:arithmetic}
For every sorted five-pile position,
\begin{equation}\label{eq:parity-facts}
 \chi=1\ \Longrightarrow\ L\text{ even},\qquad
 D\text{ even}\ \Longleftrightarrow\ Q(x).
\end{equation}
Moreover,
\begin{equation}\label{eq:regions}
 F(x)=
 \begin{cases}
 \parround(N,D),&\Delta<0,\\
 D,&\Delta\ge0,
 \end{cases}
 \qquad
 \parround(N,D)=D\quad\text{when }\Delta=0.
\end{equation}
\end{lemma}

\begin{proof}
Condition $Q(x)$ implies that $L$ is even. The assertions in
\eqref{eq:parity-facts} therefore follow directly from the
rounding rule \eqref{eq:chi}.

If $\Delta<0$, then $\parround(N,D)\le L-1\le D$. If
$\Delta=0$, then $N=2L$ and \eqref{eq:pi} gives
$\parround(N,D)=L-\chi=D$. Finally, if $\Delta>0$ and $N$ is odd,
then $\lfloor N/2\rfloor\ge L\ge D$. If $\Delta>0$ and $N$ is
even, then $N/2\ge L+1$, so $\parround(N,D)\ge L\ge D$.
\end{proof}

\subsection{Option values near the boundary}\label{subsec:five-options}

The following calculations will be used to determine the GM
exceptions. For a canonical move with indicator $\epsilon$, put
$m=\epsilon_1+\epsilon_3+\epsilon_4$. Then
\begin{equation}\label{eq:move-parameters}
        N'=N-2,\qquad L'=L-m,\qquad \Delta'=\Delta+2m-2 .
\end{equation}
Only the pair $b,e$ has $m=0$. A move with $m=2$ acts on $a,c$,
$a,d$ or $c,d$. The definition of $\Delta$ also gives
\begin{equation}\label{eq:largest-gap}
                   \Delta\ge0\quad\Longrightarrow\quad e-d\ge a+c-b .
\end{equation}

\begin{lemma}\label{lem:five-options}
If $\Delta<0$ and $N$ is even, every move from or to a $Q$-position
leads to a position of remoteness $R_5(x)-1$.
If $N$ is odd and $\Delta\le-3$, every option has remoteness
$R_5(x)-1$. When $\Delta=-1$, put $T=(N-1)/2=L-1$.
Options with $m\le1$ have value $T-1$, and those with $m=2$ have
\begin{equation}\label{eq:low-crossing}
                       R_5(x')=T-1-\chi'.
\end{equation}
If $\Delta\ge0$, every option with $\Delta'\ge0$ satisfies
\begin{equation}\label{eq:high-options}
                       R_5(x')=L-m-\chi'.
\end{equation}
Only the pair $b,e$ can cross from $\Delta\ge0$ to $\Delta'<0$,
and then
\begin{equation}\label{eq:high-crossing}
 R_5(x')=\begin{cases}L-2+\chi',&\Delta=0,\\ L-1,&\Delta=1.\end{cases}
\end{equation}
\end{lemma}

\begin{proof}
For even $N$ and $\Delta<0$ we have $\Delta'\le0$. A two-pile
move cannot preserve $Q$ at a $Q$-position. It either changes
the parity of $a$ or changes one or two of the three parities
of $b,c,d$. Hence the positions before and after any move from or
to $Q$ have opposite $D$-parities. The half-total formula in
\eqref{eq:regions} and \eqref{eq:pi-minus} give the first claim.
For odd $N$ and $\Delta\le-3$, every option remains in $\Delta'<0$
and the half-total drops by one. When $\Delta=-1$, only $m=2$
leaves this region, giving $L'=L-2$ and $D'=L-2-\chi'$.
For the last two formulas, insert \eqref{eq:move-parameters} into
\eqref{eq:regions}. A crossing requires $m=0$ and $\Delta=0$
or $1$. Its half-total bound from \eqref{eq:pi} is then
$L-2+\chi'$ or $L-1$, respectively.
\end{proof}

\Needspace{10\baselineskip}
\section{Consequences of the formula}\label{sec:consequences}

\subsection{Three and four piles}

\begin{corollary}\label{cor:small}
For $0\le u\le v\le w$, let
$N_3=u+v+w$ and $D_3=u+v-[u,v\text{ both odd}]$. Then
\begin{equation}\label{eq:three}
 R_3(u,v,w)=\min\{D_3,\parround(N_3,D_3)\}.
\end{equation}
For $0\le a\le b\le c\le d$, let
$N_4=a+b+c+d$ and $D_4=b+c-[b-a,c-a\text{ both odd}]$. Then
\begin{equation}\label{eq:four}
 R_4(a,b,c,d)=\min\{D_4,\parround(N_4,D_4)\}.
\end{equation}
\end{corollary}

\begin{proof}
Insert two initial zero piles into \eqref{eq:five} to obtain
\eqref{eq:three}. Inserting one initial zero gives
$D=b+c-[b\equiv c\not\equiv a]$, which is $D_4$.
The correspondence obtained by adding zero piles preserves all
legal moves and terminal positions, and therefore preserves
remoteness.
\end{proof}

For example, the formulas give $R_3(0,v,w)=v$ and
$R_4(a,a,a,d)=2a$. They also give $R_3(3,3,3)=4$:
here $D_3=5$ and $\parround(9,5)=4$.
The three-pile case is already covered by the general theory of
\cite{GurvichEtAl2024}; \eqref{eq:three} is a compact expression
suited to the present proof.

\subsection{Five-pile shifts}

\begin{corollary}\label{cor:five-shifts}
For every sorted five-pile position,
\begin{align}
 R_5(a,b+1,c+1,d+1,e+1)&=R_5(a,b,c,d,e)+2,
   \label{eq:five-shift-one}\\
 R_5(a,b,c,d+2,e+2)&=R_5(a,b,c,d,e)+2.
   \label{eq:five-shift-two}
\end{align}
\end{corollary}

\begin{proof}
Both shifts increase $N$ by four and $L$ by two. The first flips
the parities of $b,c,d$ together, leaving $\chi$ unchanged; the
second preserves every parity. Thus $D$ increases by two, and
\eqref{eq:pi-plus} shows that $\parround(N,D)$ also increases by two.
Apply Corollary~\ref{thm:five}.
\end{proof}

The shifts act on the indicated \emph{sorted} coordinates.
The first is a known Sprague--Grundy shift, and the second is a
known normal-play P-shift \cite[Section~2 and Theorem~3]{Chikin2021}.
The corollary specifies their exact remoteness increments.

\subsection{All outcome-preserving translations in five and six piles}
\label{subsec:all-shifts}

A \emph{P-shift} is a nonnegative nondecreasing integer vector $y$
such that $x$ is a P-position exactly when $x+y$ is a P-position,
for every sorted position $x$ of the given dimension
\cite{Chikin2021}. An \emph{additive remoteness shift} satisfies
$R_n(x+y)=R_n(x)+k$ for a fixed integer $k$ and every such $x$.

\begin{theorem}\label{thm:all-shifts}
For $n=5,6$, these two sets of shifts coincide. The increment is
$k=(\sum_i y_i)/2$ and is always even. In each dimension, the shifts
form the free commutative monoid on the four generators in the table.
\[
\begin{array}{@{}lc@{\qquad\qquad}lc@{}}
\toprule
 n=5 & k & n=6 & k\\ \midrule
 (0,1,1,1,1)&2 &(0,0,1,1,1,1)&2\\
 (0,0,0,2,2)&2 &(1,1,1,1,2,2)&4\\
 (0,0,2,2,4)&4 &(0,0,0,2,2,4)&4\\
 (2,2,2,2,4)&6 &(0,2,2,2,2,4)&6\\ \bottomrule
\end{array}
\]
The four six-pile vectors are the P-shifts listed in
\cite[Proposition~1]{Chikin2021}.
\end{theorem}

\begin{proof}
For six piles, write $S=N-2A$ and $T=N-2B$. On the slice defined by
$E(x)$ and $N\equiv3\pmod4$, both $A$ and $B$ are even. Here $h=0$ and
\eqref{eq:six} gives
\begin{equation}\label{eq:quadrant}
                x\in P\quad\Longleftrightarrow\quad S>0,\ T>0.
\end{equation}
All pairs $(S,T)\in(4\Z+3)^2$ occur on this slice. To see this,
put $g_i=x_i-x_{i-1}$, $x_0=0$, and $r=(S-T)/2$. Take
$g_1=\max(r,0)$, $g_5=\max(-r,0)$ and $g_3=g_4=0$. Choose
$g_2$ to be a nonnegative multiple of four large enough that
$g_6=S+g_2\ge0$. Then $E$ holds, $S=g_6-g_2-g_4$ and
$B-A=g_1-g_5=r$. Thus $N=S+2A\equiv3\pmod4$ as required.

If $y$ is a P-shift, its fourth multiple preserves this slice and
translates $(S,T)$ by $(4S(y),4T(y))$. A translation preserving
membership in the positive quadrant on the whole residue lattice
must vanish. Indeed, hold one coordinate sufficiently positive
before and after translation and test the other on both sides of
zero. Consequently $S(y)=T(y)=0$. The zero position is losing,
so $y$ is losing. Its total is even, and \eqref{eq:six} now gives
$E(y)$. Conversely, if $S(y)=T(y)=0$ and $E(y)$ holds, translation
by $y$ preserves $h$ and $p$, and increases $A$ and $B$ by the
same even integer $N(y)/2$. Both terms in \eqref{eq:six} increase
by this integer.

For five piles, the same argument uses one variable.
When $Q(x)$ holds and $N\equiv3\pmod4$, formula~\eqref{eq:five} gives
$x\in P$ exactly when $\Delta>0$. Every
$v\in4\Z+3$ occurs as $\Delta$ on this slice. Take
$x=(2t,2t,2t,2t,4t+v)$ with $t$ sufficiently large.
Invariance under $4y$ forces $\Delta(y)=0$. Then $N(y)$ is even,
and $y\in P$ gives $Q(y)$. Conversely, $\Delta(y)=0$ and $Q(y)$
preserve $Q(x)$ under translation and increase $L$ by the even
integer $N(y)/2$. Equation~\eqref{eq:five} gives the claimed shift.

It remains to exclude an additive shift with odd increment. In either
dimension choose the parities of all but the largest pile so that
both the initial parity vector and its translate fail $Q$ (five
piles) or $E$ (six piles). This is possible since each condition
holds on just one eighth of these vectors. Realize the chosen
parities in sorted order and take the largest pile sufficiently
large. At both positions the term controlled by the smaller piles
determines remoteness, which is odd in both cases. Thus every additive
shift has even increment and preserves the outcome.

Finally, write the gaps of $y$ as $g_i=y_i-y_{i-1}$, $y_0=0$.
For six piles, the conditions just proved are
\[
 g_5=g_1,\qquad g_6=g_2+g_4,\qquad g_2,g_4\ \text{even}.
\]
The parameters $g_1,g_2/2,g_3,g_4/2$ range independently over the
nonnegative integers. For five piles, the corresponding conditions are
\[
 g_5=g_1+g_3,\qquad g_1,g_3,g_4\ \text{even}.
\]
Now $g_1/2,g_2,g_3/2,g_4/2$ are free. The four unit vectors in each
parameter space give the displayed generators. Each shift has a unique
representation in terms of these generators.
\end{proof}

\subsection{Optimal moves and evaluation}

Corollary~\ref{thm:five} determines every remoteness-optimal move:
evaluate the formula at each of the at most ten legal successors
and retain those with value $R_5(x)-1$. It also identifies
P-positions by the parity of its value. Explicitly,
\begin{equation}\label{eq:outcome}
 x\text{ is a P-position}\quad\Longleftrightarrow\quad
 \begin{cases}
 Q(x),&N\text{ even or }\Delta\ge0,\\
 N\equiv1\pmod4,&N\text{ odd and }\Delta<0.
 \end{cases}
\end{equation}
This recovers the outcome description in
\cite[Corollaries~1--2]{Chikin2021}; the remoteness formula
refines it by giving the duration of optimal play.
For fixed dimension,
the formula uses a bounded number of additions, comparisons,
parity tests and halvings. Hence it can be evaluated in $O(B)$
bit operations when each pile size has at most $B$ bits.
This complexity statement concerns the arithmetic formula and
does not require enumerating the game tree.
Corollary~\ref{cor:corrected-gm} will select one optimal move
in the same time bound without evaluating any successor.

\section{The GM rule and its exceptions}\label{sec:gm}

For $n\ge2$, we use the following deterministic GM rule. Keep
$n-2$ piles unchanged and decrease the remaining two. If at least
$n-2$ piles are even, keep the $n-2$ smallest even piles. Otherwise,
keep all even piles and fill the remaining places with the largest
odd piles. Among equal piles, keep those with the rightmost indices.
This convention makes the move canonical. It is legal at every
nonterminal position because all zero piles are kept. It is the
largest-remaining-pile choice allowed in
\cite[Section~2.1]{GurvichMaximchukNaumova2026}; related move rules
were developed in \cite{GurvichNaumova2023,GurvichNaumova2024}.

We call a position \emph{exceptional} when its GM move is not
remoteness-optimal. At a winning position, such a move may give up
the win or delay it. At a losing position, it shortens the time to defeat.

\Needspace{12\baselineskip}
\subsection{Five piles with at least three even piles}

\begin{theorem}\label{thm:gm-five}
Let $x=(a,b,c,d,e)$ be a sorted nonterminal position with at least
three even piles. It is exceptional if and only if
\begin{equation}\label{eq:five-exception}
 a,c,d\text{ are even},\qquad b\text{ is odd},\qquad
 e\ge a+c+d-b+1.
\end{equation}
At an exception,
\begin{equation}\label{eq:exception-value}
 r:=R_5(x)=a+c+d-1,
\end{equation}
and the unique remoteness-optimal move decreases $c,d$.
The GM move decreases $b,e$. Every legal successor has the
remoteness listed in Table~\ref{tab:exceptional-successors}.
\end{theorem}

\begin{table}[htbp]
\centering
\small
\renewcommand{\arraystretch}{1.2}
\begin{tabular}{@{}p{0.76\textwidth}c@{}}
\toprule
Piles decreased and, when needed, the boundary condition & Remoteness\\
\midrule
$c,d$ & $r-1$\\
$b,e$, with $e=a+c+d-b+1$ & $r$\\
$b,e$, with $e>a+c+d-b+1$ & $r+1$\\
$a,c$ or $a,d$ & $r-2$\\
$a,b$; $a,e$; $b,c$; $b,d$; $c,e$; or $d,e$ & $r$\\
\bottomrule
\end{tabular}
\caption{Successor remotenesses under \eqref{eq:five-exception}.
Moves involving $a$ are omitted when $a=0$.}
\label{tab:exceptional-successors}
\end{table}

\begin{proof}
Let $q$ be the number of odd piles. By assumption $q\in\{0,1,2\}$.
We first determine precisely when the GM move fails to reduce
remoteness by one.

\smallskip
\noindent\emph{No odd piles.}
If $\Delta<0$, Lemma~\ref{lem:five-options}
shows that every successor has remoteness $R_5(x)-1$.
Suppose $\Delta\ge0$. Then $Q(x)$ holds and $R_5(x)=L$.
The GM move decreases two largest piles. Its canonical
representative has $m\in\{0,1\}$: the alternative $m=2$ would
require $b<c=d=e$, giving $\Delta=b-a-c<0$. If all four largest
piles are equal, $\Delta\ge0$ forces $a=0$, and the canonical
pair is $b,c$, with $m=1$.
For $m=0$ we have $\chi'=1$, and for $m=1$ we have $\chi'=0$.
Equations~\eqref{eq:high-options} and \eqref{eq:high-crossing}
therefore give successor remoteness $L-1$ in every case.

\smallskip
\noindent\emph{Two odd piles.}
The GM move decreases those two piles and gives an all-even
successor. If $\Delta<0$, Lemma~\ref{lem:five-options} gives remoteness
$R_5(x)-1$. Suppose $\Delta\ge0$, and write $r=R_5(x)$.
The original position fails $Q$, so $r=D$ is odd, whereas
$\chi'=0$. If $L$ is odd, making $L'$ even requires $m=1$,
and $D'=r-1$. If $L$ is even, $\chi=1$ and $m$ is zero or
two. The choice $m=2$ gives $D'=r-1$. The choice $m=0$
occurs exactly when $b,e$ are the odd piles. It gives
remoteness $r-1$ at $\Delta=0$, by
\eqref{eq:high-crossing}, and $r+1$ at $\Delta\ge2$.
These are precisely the positions satisfying \eqref{eq:five-exception}
with $e$ odd.

\smallskip
\noindent\emph{One odd pile and $\Delta<0$.}
The GM move decreases the odd pile and a largest even pile.
If $\Delta\le-3$, every successor has remoteness $R_5(x)-1$.
At $\Delta=-1$, only an $m=2$ move could give a different value,
as in \eqref{eq:low-crossing}. Such a canonical move excludes
the indices of $b,e$. Consequently the largest even pile
selected by the rule must be tied with $e$, and $d=e$ is even.
The equality $\Delta=-1$ becomes
\[
 a+(c-b)=1.
\]
Either $a=1$ and $c=b$, or $a=0$ and $c=b+1$.
In the first case $a$ is the unique odd pile. In the second,
an $m=2$ move forces $c$ to be the unique odd pile, since
the odd pile must be selected and neither $b$ nor $e$ is
selected. Thus $L$ is odd in both cases. It follows that
$L'=L-2$ is odd and $\chi'=0$, so
\eqref{eq:low-crossing} still gives $R_5(x)-1$.

\smallskip
\noindent\emph{One odd pile and $\Delta\ge0$.}
If the odd pile is $b$, then $a<b<c$ and $e>d$, the latter by
\eqref{eq:largest-gap}. The GM move decreases $b,e$.
Here $\chi=1$ and $\Delta$ is positive and odd. The GM successor
has remoteness $r$ at $\Delta=1$, and $r+1$ at
$\Delta\ge3$. These are precisely the exceptions in
\eqref{eq:five-exception} with $e$ even.

If the odd pile is $a$, $c$ or $d$, the GM move decreases that
pile and $e$. In the first two cases $e>d$ by
\eqref{eq:largest-gap}; in the third this follows from parity.
The move has $m=1$ and leads to a $Q$-position of remoteness $r-1$.
If the odd pile is $e$, then $Q(x)$ holds and $r=L$.
The GM move decreases $e$ and a largest even pile. Its
canonical representative has $m=0$ or $1$, and the same
calculation as in the all-even case gives remoteness $r-1$.
This proves the equivalence in \eqref{eq:five-exception}.

It remains to compute all successor values at an exception.
There $\chi=1$, $\Delta\ge1$, and $r=L-1$, proving
\eqref{eq:exception-value}. The ordering is
\[
 a<b<c\le d<e.
\]
If $c=d$, a move selecting just one of them may be represented
by selecting $c$; the resulting sorted position is the same.
Decreasing $c,d$ gives $m=2$ and $\chi'=0$, hence value $r-1$.
Decreasing $a,c$ or $a,d$, when legal, gives $m=2$ and
$\chi'=1$, hence value $r-2$. Each of the six listed $m=1$
moves gives odd $L'=L-1$ and $\chi'=0$, hence value $r$.
The $m=0$ move on $b,e$ gives $r$ when $\Delta=1$, by
\eqref{eq:high-crossing}, and otherwise gives $r+1$.
This proves the table and uniqueness of the optimal move.
\end{proof}

The family \eqref{eq:five-exception}, its remoteness, the unique
optimal move, and the two values of the GM successor were
proposed in \cite[Section~2.3]{GurvichMaximchukNaumova2026}.
The theorem proves those statements and supplies the remaining
successor values.

\begin{example}\label{ex:boundary}
For $x=(2,3,4,4,8)$, the example after Corollary~\ref{thm:five}
gives $R_5(x)=9$ and $\Delta=1$. Decreasing the two piles of size
four leaves a position of remoteness eight. The GM successor has remoteness
nine and is an N-position, so the GM move loses.
Increasing the largest pile from eight to nine gives $\Delta=2$ and
leaves the remoteness at nine. The GM successor now has remoteness ten.
With optimal play thereafter, the GM move gives a win in eleven moves
in total. An optimal first move gives a win in nine.
\end{example}

\subsection{The complete classification and an optimal rule}

The positions with fewer than three even piles introduce a second
exception family. Keeping the largest odd piles forces the GM move
to decrease the two smallest odd piles.

\begin{theorem}[All five-pile exceptions]\label{thm:gm-complete}
A sorted nonterminal position $x=(a,b,c,d,e)$ is exceptional for
the specified GM rule if and only if it belongs to one of the
following two disjoint families:
\begin{enumerate}
\item[(I)] $a,c,d$ are even, $b$ is odd, and $\Delta>0$;
\item[(II)] $a$ is odd, at least two of $b,c,d$ are odd, and
      either $N$ is even or $\Delta>0$.
\end{enumerate}
Every exception is an N-position. In family~\textup{(II)} the GM
successor is also an N-position, so the GM move loses. In
particular, the GM rule is remoteness-optimal at every P-position.
\end{theorem}

\begin{proof}
Theorem~\ref{thm:gm-five} handles positions with at least three
even piles, and its exception condition is exactly family~(I).
Suppose now that at least three piles are odd. The GM move decreases
the first two odd coordinates; this is a legal canonical move.
Let $x'$ be its successor and put
\[
 H:\quad a\text{ is odd and at least two of }b,c,d\text{ are odd}.
\]
The change in $D$ is
\begin{equation}\label{eq:gm-low-even-change}
 D'=\begin{cases}
 D-2,&a,c,d\text{ odd},\\
 D,&a,b\text{ odd and }c\not\equiv d,\\
 D-1,&H\text{ fails}.
 \end{cases}
\end{equation}
The first two cases partition the positions satisfying $H$.
In the first case, $c,d$ are odd, so $D=L$.
The move decreases $a,b$ when $b$ is odd,
or $a,c$ when $b$ is even. The respective values of
$(L',\chi')$ are $(L-1,1)$ and $(L-2,0)$.
If $c,d$ have opposite parities and $a,b$ are odd, then
$D=L-1$ and decreasing $a,b$ gives $(L',\chi')=(L-1,0)$.

For the last line, first suppose $a$ is even. If $b,c,d$ are
all odd, the move on $b,c$ changes $(L,\chi)$ to $(L-1,0)$,
with $\chi=0$. Otherwise exactly two of $b,c,d$ are odd and
$e$ is odd. If those two are $c,d$, the move changes
$(L,1)$ to $(L-2,0)$; if one is $b$, it changes
$(L,0)$ to $(L-1,0)$. Finally, if $a$ is odd and $H$ fails,
exactly one of $b,c,d$ is odd, as is $e$. When that coordinate
is $b$, the change is $(L,0)\mapsto(L-1,0)$; otherwise it is
$(L,1)\mapsto(L-2,0)$. This proves
\eqref{eq:gm-low-even-change} in all cases.

If $H$ fails, the identities $N'=N-2$, $D'=D-1$ and
\eqref{eq:pi-minus} give $R_5(x')=R_5(x)-1$ immediately.
Suppose $H$ holds. Then $D,D'$ are both odd. If $N$ is even,
the formula shows that both $R_5(x)$ and $R_5(x')$ are odd,
so the move is exceptional and loses.

If $N$ is odd, put $T=(N-1)/2$. The two remotenesses are
\[
 R_5(x)=\min(D,T),\qquad R_5(x')=\min(D',T-1).
\]
For $\Delta<0$ we have $T\le L-1$. When $D=L$ and
$D'=L-2$, these minima are $T,T-1$; when $D=D'=L-1$,
they are again $T,T-1$. Thus the move is optimal.
For $\Delta>0$ we have $T\ge L$. The first remoteness is
$D$, while the second is $D-2$ or $D$, according to the
first two lines of \eqref{eq:gm-low-even-change}. Both are
odd, so the move loses. Since odd $N$ excludes $\Delta=0$,
these are precisely the positions in family~(II).
Family~(I) also has odd remoteness by
\eqref{eq:exception-value}. This proves all assertions.
\end{proof}

\begin{corollary}[A corrected GM rule]\label{cor:corrected-gm}
The following rule is remoteness-optimal at every nonterminal
sorted five-pile position:
\begin{enumerate}
\item In family~\textup{(I)}, decrease $c,d$.
\item In family~\textup{(II)}, decrease $a$ and a smallest even pile.
      If all piles are odd, decrease $a,e$ instead.
\item At every other position, use the GM move.
\end{enumerate}
Among equal piles use the canonical representative.
\end{corollary}

\begin{proof}
Theorem~\ref{thm:gm-five} proves the first instruction, and
Theorem~\ref{thm:gm-complete} proves the third.
In family~(II), $a$ is positive and every even pile is strictly
larger than $a$, so the second instruction is legal.
Suppose first that $\Delta\ge0$. After the corrected move, $a$ is
even and $b,c,d$ have the same parity. If $b$ is the even member of $b,c,d$,
or if all three are odd, $L$ is odd and the pair is respectively
$a,b$ or $a,e$. Thus $L'=L-1$. When $c$ or $d$ is the even member,
$L$ is even, $\chi=1$, and the correction gives $L'=L-2$.
In each case $D'=D-1$ and $m\ge1$ ensures $\Delta'\ge0$.
A selected even pile among $b,c,d$ starts its block.
When $e$ is selected, \eqref{eq:largest-gap} gives $e>d$
because $a>0$. The moves are therefore canonical.

If $\Delta<0$, family~(II) requires even $N$. There are
exactly four odd piles, so the rule decreases $a$ and the
unique even pile. The successor satisfies $Q$. Since
$\Delta\le-2$, its imbalance is at most zero by
\eqref{eq:move-parameters}. The even-total calculation in
Lemma~\ref{lem:five-options} gives its
remoteness as $R_5(x)-1$.
\end{proof}

For example, $(1,1,1,1,2)$ has remoteness three and belongs
to family~(II). The GM move reaches $(0,0,1,1,2)$, of
remoteness one, whereas the corrected move reaches
$(0,1,1,1,1)$, of remoteness two. Thus family~(II) is needed
even for small positions.

\begin{corollary}[Allowing all choices of odd pivots]\label{cor:all-gm}
Allow the extra odd piles kept by the GM rule to be chosen arbitrarily,
as in \cite[Section~2.1]{GurvichMaximchukNaumova2026}. An optimal
GM move exists at every five-pile nonterminal position except those
in family~\textup{(I)} and those in family~\textup{(II)} containing
an even pile.
\end{corollary}
\begin{proof}
Outside the two families the deterministic choice already works. In
family~(I) the rule fixes the move up to equal piles. When all piles
are odd, the correction $a,e$ in family~(II) is an admissible choice.
It remains to consider family~(II) with an even pile. Every admissible
GM move decreases two odd piles. If the total is even, there are four
odd piles before the move and two after it. No even-total P-position
has two odd piles, by \eqref{eq:outcome}.

If the total is odd, there are three odd piles, namely $a$ and two
of $b,c,d$. Here $e$ is even and $\Delta>0$. A move on two odd
piles has $m\ge1$ in its canonical representative, so $\Delta'>0$.
If $a$ is left odd, the successor fails $Q$. If $a$ is decreased,
exactly one of the two odd entries among $b,c,d$ is also decreased
and the three entries still include both parities. In either case,
the successor is winning. Since the original position is also winning,
an optimal move would lead to a losing position.
\end{proof}

\subsection{A complete four-pile classification}

\begin{theorem}\label{thm:gm-four}
A sorted nonterminal four-pile position $(a,b,c,d)$ is
exceptional for the specified GM rule if and only if
\begin{equation}\label{eq:four-exception}
 a\text{ is odd},\qquad b,c\text{ are even},\qquad
 d\ge b+c-a+1.
\end{equation}
At an exception, $r:=R_4(a,b,c,d)=b+c-1$, and the unique
remoteness-optimal move decreases $b,c$. The GM move decreases
$a,d$; its successor has remoteness $r$ when
$d=b+c-a+1$ and $r+1$ when $d>b+c-a+1$.
Each of the other four choices of a pair leads to remoteness $r$.
When $b=c$, distinct choices can give the same sorted successor.
\end{theorem}

\begin{proof}
Insert a zero pile and apply Theorem~\ref{thm:gm-complete} to
$(0,a,b,c,d)$. The five-pile GM move keeps the zero and
exactly the two piles kept by the four-pile rule, up to the
interchange of equal zero piles. Family~(II) is impossible,
and family~(I) becomes \eqref{eq:four-exception}.
Table~\ref{tab:exceptional-successors} gives all successor
values: moves involving the inserted zero are illegal, and
the six remaining pairs give precisely the stated values.
This also proves the pattern proposed in
\cite[Section~2.2]{GurvichMaximchukNaumova2026}.
\end{proof}

For instance, $(1,2,2,4)$ is an exception of remoteness three.
The optimal move decreases the middle piles, whereas the GM
move loses. By contrast, $(1,1,1,2)$ has remoteness two and
an optimal GM move, despite having only one even pile.

\subsection{When a six-pile move ceases to be optimal}
\label{subsec:six-gm}

Once the five smaller piles have been fixed, the largest pile enters
\eqref{eq:six} only through the half-total bound. This gives a
threshold for the optimality of each move, including the GM move.

\begin{theorem}[A threshold for each pair]\label{thm:pair-threshold}
Let a legal canonical pair have row $v$ and losses $(m,n)$ in
Table~\ref{tab:six-pairs}. Put $u'=u\oplus v$, $r=B-A$, and
\begin{align}
 K={}&1+\max\{-m-\alpha(u'),\ r-n-\beta(u')\}
              -\max\{-\alpha(u),\ r-\beta(u)\},
 \label{eq:pair-excess}\\
 s={}&\parround(N,D_6)-D_6.\notag
\end{align}
If $N$ is even and $K$ is odd, the pair fails to decrease
remoteness by one.
In every other case it is optimal exactly when
\begin{equation}\label{eq:cutoff}
                 K=0\qquad\text{or}\qquad s\le\min(0,K).
\end{equation}
Fix the first five piles, write $S=a+b+c+d+e$, and let $f>e$ vary
through one parity class. For $K\ne0$, except when $N$ is even and
$K$ is odd, the move is suboptimal precisely when
\begin{equation}\label{eq:pair-threshold-f}
       f\ge 2\bigl(D_6+\min(0,K)+1+[N\text{ even}]\bigr)-S .
\end{equation}
This also gives the cutoff for the deterministic GM pair, which is
fixed on each such class of nonterminal positions.
\end{theorem}

\begin{proof}
Subtract $A$ from the two maxima in \eqref{eq:offset-form} before
and after the move. This gives $K=D'_6-D_6+1$. If $N$ is even
and $K$ is odd, $D'_6$ and $D_6$ have the same parity. Their
remotenesses also have that parity, so the move cannot be optimal.
Otherwise $\parround(N-2,D'_6)=\parround(N,D_6)-1$.
For odd $N$ this follows from the half-total, and for even $N$
it follows from the opposite parities and \eqref{eq:pi-minus}.
Taking the two minima in \eqref{eq:six} now gives the exact error
\[
             R_6(x')-R_6(x)+1=\min(K,s)-\min(0,s).
\]
It vanishes precisely under \eqref{eq:cutoff}.
For odd $N$, the inequality $s>\min(0,K)$ requires
$\lfloor N/2\rfloor\ge D_6+\min(0,K)+1$. For even $N$ and even
$K$, both $s$ and $\min(0,K)$ are even, so the required bound on
$N/2$ is $D_6+\min(0,K)+2$. This gives \eqref{eq:pair-threshold-f}.
With $f>e$ and its parity fixed, the GM choice has constant indices
and hence constant $K$.
\end{proof}

\begin{corollary}[All-even exceptions]\label{cor:even-six}
A sorted nonterminal six-pile position in which every pile is even
is exceptional for the GM rule if and only if
\begin{equation}\label{eq:even-six-exception}
              a>0,\qquad c<d=e,\qquad f\ge a+b+2d-c .
\end{equation}
Every such position is losing. Its GM successor has remoteness
$R_6(x)-3$.
\end{corollary}

\begin{proof}
Here $u=000$ and both $A$ and $B$ are even. Suppose first that $f>e$.
The GM rule decreases pile~6 and the first pile of size $e$.
By Table~\ref{tab:six-pairs}, each such pair lowers $X$
by one. It lowers $Y$ by one as well except for $46$, which lowers
$Y$ by three. This pair is chosen exactly when $c<d=e$. If $B\le A$,
the maximum still falls by one. If $B>A$, equivalently $a>0$
when $d=e$, then $B\ge A+2$, $D_6=B$ and $K=-2$.
Since $s=\rho_0(\lfloor N/2\rfloor)-B$ is even, the move is
exceptional precisely when $s\ge0$, or $N\ge2B$. Expanding $N$
and using $e=d$ yields \eqref{eq:even-six-exception}.
The error formula in the theorem then gives $R_6(x')=R_6(x)-3$.

If $f=e$, the all-even position is a P-position with tied largest
piles. Every move is optimal by Corollary~\ref{cor:tied}.
\end{proof}

\begin{example}
The position $(2,2,2,4,4,10)$ has remoteness $12$.
Decreasing a pile of size $4$ and the pile of size $10$ is the
GM move and leaves remoteness $9$. Decreasing a pile of size $2$
and the pile of size $10$ leaves remoteness $11$. The losing
player can therefore prolong play by two moves.
\end{example}

\section{Complementary moves and remoteness shifts}\label{sec:shifts}

The shift identities in this section do not require the five-pile
formula. In each proof, a move from a shifted position either
corresponds to a move from the original position or allows a reply
that returns to the original position.

\begin{lemma}[Two-step shift criterion]\label{lem:two-step}
Let $\Phi$ be a map on the positions of a finite-play impartial
game as in Section~\ref{sec:preliminaries}. Make the following assumptions.
\begin{enumerate}
\item $x\to y$ implies $\Phi(x)\to\Phi(y)$;
\item for every move $\Phi(x)\to z$, either
      $z=\Phi(y)$ for some move $x\to y$, or $z\to x$;
\item if $x$ is terminal, then $\Phi(x)$ is nonterminal.
\end{enumerate}
Then $R(\Phi(x))=R(x)+2$ for every position $x$.
\end{lemma}

\begin{proof}
Induct on the maximum length of play from $x$, and put $r=R(x)$.
An option of the form $\Phi(y)$ has remoteness $R(y)+2$ by
induction. From any other option, a move back to $x$ is legal.

If $r=0$, then $x$ is terminal. Every option of $\Phi(x)$ has
remoteness one because it can move to $x$. The third hypothesis
therefore gives $R(\Phi(x))=2$.

If $r>0$ is even, every original option $y$ has odd remoteness
at most $r-1$. Hence every translated option has odd remoteness
at most $r+1$. An option from which a move to $x$ is legal also has odd
remoteness at most $r+1$, by the minimum clause of
\eqref{eq:recurrence}. Translating an optimal move from $x$
attains $r+1$, so the maximum clause gives $R(\Phi(x))=r+2$.

If $r$ is odd, every even-remoteness translated option has value
at least $r+1$. The same bound holds for an even-remoteness
option that can move to $x$: the maximum clause of the recurrence
gives a value at least $1+R(x)$. A translated optimal move
attains the even value $r+1$. The minimum clause therefore
gives $R(\Phi(x))=r+2$.
\end{proof}

\begin{theorem}[Complement shift]\label{thm:complement}
Let $S\subseteq\Z_{\ge0}^{d}\setminus\{0\}$ be nonempty,
and suppose that, for some $c\in\Z_{\ge0}^{d}$,
\begin{equation}\label{eq:complement-condition}
 S=c-S:=\{c-v:v\in S\}.
\end{equation}
In the vector subtraction game with legal moves
$x\to x-v$ for $v\in S$, $v\le x$, one has
\begin{equation}\label{eq:complement-shift}
 R(x+c)=R(x)+2\qquad(x\in\Z_{\ge0}^{d}).
\end{equation}
\end{theorem}

\begin{proof}
Condition~\eqref{eq:complement-condition} implies $v\le c$
for every $v\in S$. Hence $S$ is finite, every move strictly
decreases the coordinate sum, and $x+c$ is always nonterminal.
Translation by $c$ preserves every legal move. Moreover, every
option $z=x+c-v$ has the legal complementary reply
\begin{equation}\label{eq:reply}
 z\longrightarrow z-(c-v)=x.
\end{equation}
All hypotheses of Lemma~\ref{lem:two-step} hold for
$\Phi(x)=x+c$.
\end{proof}

For a two-element subtraction set $S=\{u,v\}$, take $c=u+v$.
The theorem then refines the two-move outcome translation
principle of Larsson, Saha and Yokoo
\cite[Lemma~2]{LarssonSahaYokoo2023} by giving its exact
remoteness increment.

\begin{corollary}[Complementary cardinalities]\label{cor:cardinalities}
Let $A\subseteq\{1,\ldots,n-1\}$ be nonempty. Suppose a move
removes one token from each of exactly $j$ distinct positive
piles, for any $j\in A$. If
\[
 A=n-A:=\{n-j:j\in A\},
\]
then $R(x+\one)=R(x)+2$.
\end{corollary}

\begin{proof}
The complement of an allowed subset of piles is again allowed.
Thus the set $S$ of allowed $0$--$1$ subtraction vectors satisfies $S=\one-S$,
and Theorem~\ref{thm:complement} applies.
\end{proof}

For example, the corollary applies to the five-pile game in
which one removes a single token from each of either two or
three piles. The restriction to unit removal is part of the
rules. For exact slow $k$-Nim with $2k$ piles, it gives
\begin{equation}\label{eq:diagonal}
 R_{2k,k}(x+h\one)=R_{2k,k}(x)+2h\qquad(h\in\Z_{\ge0}),
\end{equation}
where $R_{n,k}$ denotes remoteness when exactly $k$ piles
are decreased. In particular,
\begin{equation}\label{eq:four-reduction}
 R_4(a,b,c,d)=2a+R_3(b-a,c-a,d-a).
\end{equation}
Chikin et al.~\cite[Theorem~3]{Chikin2021} proved the corresponding
Sprague--Grundy invariance.

\begin{proposition}[Two further shift families]\label{prop:tail-shifts}
For sorted positions and $k\ge1$,
\begin{align}
 R_{2k+1,k}\bigl(x+(0,1^{2k})\bigr)&=R_{2k+1,k}(x)+2,
    \label{eq:tail-one}\\
 R_{k+1,k}\bigl(x+(0,2^k)\bigr)&=R_{k+1,k}(x)+2.
    \label{eq:tail-two}
\end{align}
Here $1^{2k}$ and $2^k$ denote repeated coordinates.
\end{proposition}

\begin{proof}
Let $\Phi$ be either translation in the statement, and work
with sorted positions. Every canonical move $x\to t$ gives
a move $\Phi(x)\to\Phi(t)$. Indeed, the first coordinate of
$t$ is still a minimum by Lemma~\ref{lem:canonical}; adding
the same positive amount to all other coordinates preserves
their order and keeps the first coordinate smallest.

Conversely, consider an option of $\Phi(x)$. If its move
decreases the first coordinate, then $x_1>0$, so the same
move is legal at $x$. The first coordinate remains smallest,
and the resulting option is $\Phi(t)$ for an original
successor $t$, after sorting the last coordinates if needed.
If the move does not decrease the first coordinate, it has
a legal reply back to $x$. For \eqref{eq:tail-one}, decrease
the complementary $k$-subset of the last $2k$ coordinates.
For \eqref{eq:tail-two}, decrease each of the last $k$ coordinates
once more.

Finally, $\Phi(x)$ has at least $k$ positive coordinates in
both games and is therefore nonterminal. Apply
Lemma~\ref{lem:two-step}.
\end{proof}

These two identities refine the corresponding
Sprague--Grundy shifts of \cite[Theorem~3]{Chikin2021}.
For $k=2$, \eqref{eq:tail-one} also gives a second proof of
\eqref{eq:five-shift-one}.

\subsection{Arbitrarily many unit piles}\label{subsec:unit-piles}

The complement shift can also change the number of piles. Consider
$s$ unit piles together with a list $y$ of distinguished piles.
If fewer than $k$ piles are present, the position is terminal.

\begin{proposition}\label{prop:unit-shift}
For $k\ge2$, $1\le d\le k$, $y\in\mathbb Z_{\ge0}^d$ and $s\ge0$,
\[
          R_{s+2k,k}(1^{s+2k-d},y+\one)=R_{s+d,k}(1^s,y)+2 .
\]
\end{proposition}
\begin{proof}
Record a position by $(s,y)$. Choosing the distinguished indices
$J\subseteq\{1,\ldots,d\}$ removes the vector
$(k-|J|,\one_J)$, where $\one_J$ is its incidence vector.
A move is legal exactly when this vector is at most $(s,y)$
coordinatewise. The vectors for $J$ and its complement sum to
$(2k-d,\one)$, so Theorem~\ref{thm:complement} applies.
\end{proof}

\begin{corollary}[At most two piles larger than one]\label{cor:two-large}
In exact slow $2$-Nim put
\[
 L(s,w)=\min\left\{\left\lfloor\frac{s+w}{2}\right\rfloor, 2\left\lfloor\frac{s}{3}\right\rfloor+[3\nmid s]\right\}.
\]
Then $R_{s+1}(1^s,w)=L(s,w)$. For two distinguished piles $0\le u\le f$, put
$t=\min(\lfloor s/2\rfloor,u)$ and
$(r,v,w)=(s-2t,u-t,f-t)$. For every $s$,
\[
 R_{s+2}(1^s,u,f)=2t+
 \begin{cases}
 L(r,w),&v=0,\\
 v+[r=1,\ v\text{ even},\ w>v],&v>0 .
 \end{cases}
\]
\end{corollary}
\begin{proof}
For one distinguished pile the shift is $(3,1)$, and each term
in the minimum defining $L$ increases by two under this shift.
If $w=0$, play lasts $\lfloor s/2\rfloor$ moves, which is the
smaller bound (write $s=3q+j$ with $0\le j\le2$).
If $s\le2$ and $w>0$, the value is zero for $s=0$ and one otherwise.
These boundary values agree with $L$ and determine it everywhere.

For two distinguished piles, apply the shift $(2,1,1)$ a total
of $t$ times in reverse. When $v=0$, the one-pile formula applies.
Otherwise $r=0$ or $1$. When $r=0$, the only possible play lasts $v$ moves.
With $r=1$, formula~\eqref{eq:three} at $(1,v,w)$ gives $v$ if
$w=v$, and $1+v-[v\text{ odd}]$ if $w>v$. These are exactly
the two possibilities in the second line.
\end{proof}

\section{Further questions}\label{sec:questions}

For exact slow $2$-Nim on seven piles, Theorem~\ref{thm:transfer} separates two
questions. Can six-pile Moore remoteness be evaluated explicitly?
In the even-total region $2D>N$, which positions attain $N/2$
rather than $N/2-1$? The example in Section~\ref{sec:six} shows
why ordered pile parities alone cannot characterize the even-total
outcome for seven piles. We also leave open whether every global
P-shift for seven piles adds a constant to remoteness.
 
\section*{Declaration of AI assistance}
The authors used ChatGPT to assist with symbolic simplification of the parity
and threshold formulas, proof refinement, exposition, and the accompanying
verification scripts. Both authors independently checked and verified all
proofs and results.

\section*{Supplementary material}
All mathematical results above are proved in the text. The supplementary
programs independently test the general reduction, the explicit formulas
and the optimal move rules against Smith's recurrence. The calculations
cover all 906,192 sorted six-pile positions in which the largest pile
has size at most~26. The pair criterion agrees with the recurrence for
all 9,619,155 canonical moves from these positions.

\ifdefined\ShowFiniteCheck
  \appendix
  \input{finite_check}
\fi

\Needspace{9\baselineskip}

\end{document}

%% file: finite_check.tex
\section{An independent finite verification}\label{sec:six-proof}

Let $F$ denote the right-hand side of \eqref{eq:six}. We verify
the recurrence \eqref{eq:recurrence} for $F$. At a terminal position
$(0,0,0,0,0,f)$, the formula gives zero. At every nonterminal position,
$e>0$, so $M\ge A\ge1$ and $N\ge2$. Hence $F\ge0$. A
\emph{local check} compares $F(x)$ with one plus the minimum even
option value, or with one plus the maximum option value if all are
odd. The check at a nonterminal position is invariant under adding
the same even integer to the value at $x$ and to every option value.

\subsection*{Gaps and the legal pairs}
Put $g_i=x_i-x_{i-1}$ with $x_0=0$. A canonical pair with indicator
$\epsilon$ is legal if and only if
\[
             g_i+\epsilon_{i-1}-\epsilon_i\ge0
            \quad(1\le i\le6),\qquad \epsilon_0=0.
\]
Since these are unit moves, the list of legal pairs depends only on
which gaps vanish. For each gap, record $\ell_i=0$ if it vanishes,
$\ell_i=1$ if it is odd, and $\ell_i=2$ if it is positive and even.
A type $\ell\in\{0,1,2\}^6$ fixes the legal pairs and all six pile parities.
Write
\[
 r=B-A=g_1-g_5,\qquad s=N-2A=g_6-g_2-g_4.
\]
Within a type, $r$ and $s$ vary independently in steps of two over
their respective intervals. Indeed, $g_i=\ell_i+2t_i$, where $t_i=0$ for $\ell_i=0$ and
$t_i$ is any nonnegative integer otherwise. Let
$r_0=\ell_1-\ell_5$ and $s_0=\ell_6-\ell_2-\ell_4$. The lower and
upper endpoints of the $r$ interval are $r_0$ except that the lower
endpoint is $-\infty$ if $\ell_5>0$ and the upper is $+\infty$ if
$\ell_1>0$. Similarly, the endpoints for $s$ are $s_0$, with lower
endpoint $-\infty$ if $\ell_2+\ell_4>0$ and upper endpoint $+\infty$
if $\ell_6>0$. These descriptions include singleton intervals.

\subsection*{Which values can affect a local check}
In a move, $r$ changes by $-\epsilon_1-\epsilon_4+\epsilon_5$, which
lies between $-2$ and $1$. Let $\mathcal Z$ consist of the parent position
and its canonical options. If $r\le-1$ at the parent, every $z\in\mathcal Z$
satisfies $M(z)=A(z)$. If $r\ge3$ at the parent, every such $z$ satisfies
$M(z)=B(z)-h(z)$. Split the parent positions into five regions
\[
             r\le-1,\quad r=0,\quad r=1,\quad r=2,\quad r\ge3.
\]
Take $J=B(x)$ in the last region and $J=A(x)$ in the other four.
Within a type and a region, $J\bmod2$ and all indexed values
$D_6(z)-J$ are fixed. For the three singleton regions, this follows
by writing $M(z)=A(z)+\max(0,r(z)-h(z))$. In each outer region,
one linear expression applies at every $z$. A move decreases
$A$ and $B$ by
$\epsilon_2+\epsilon_4+\epsilon_5$ and
$\epsilon_1+\epsilon_2+2\epsilon_4$, respectively. The parity adjustment of each
option is fixed by the type. These observations also give the uniform
bounds
\begin{equation}\label{eq:offset-bounds}
                    -5\le D_6(z)-J\le3\qquad(z\in\mathcal Z).
\end{equation}

Set $q=N(x)-2J$. Let $j=0$ for the parent and $j=1$ for an option.
Since $N(z)=N(x)-2j$, the other term in $F(z)$, after subtracting $J$,
is $\lfloor(q-2j)/2\rfloor$ with a possible subtraction of one.
Its parity adjustment depends on $q\bmod4$, the parity of $J$ and
the predetermined parity $p(z)$. For $q\ge20$, inequality
\eqref{eq:offset-bounds} makes $D_6$ the minimum at every $z$. For
$q\le-20$, the half-total term is the minimum at every $z$.
Changing $q$ by a multiple of four within this lower range adds the
same even integer to all the normalized values. Hence one
representative of each attainable residue modulo four suffices in
each tail. In the middle, it suffices to check $-19\le q\le19$.

\subsection*{The finite enumeration}
For each type, intersect the interval for $r$ with the region,
preserving $r\equiv r_0\pmod2$. Skip empty intersections and denote the
endpoints of each nonempty intersection by $r_-,r_+$. For a residue $v$ of $q$ modulo four,
take the part $[s_-,s_+]$ of the $s$ interval with
$s\equiv v\pmod4$ if $J=A$, or $s\equiv v+2r_0\pmod4$ if $J=B$.
Again, skip empty intersections. The attainable values of $q$ with
this residue form the step-four interval
\[
 [s_-,s_+] \quad(J=A),\qquad
 [s_--2r_+,s_+-2r_-]\quad(J=B).
\]
The second interval follows from $q=s-2r$ and the independence of
$r$ and $s$. The sum of two intervals with step size four contains
every point of the resulting residue class between its endpoints.

For clarity, one may evaluate the local expressions without
reconstructing large reachable pile sizes. Choose any representative
$r_*$ of the region in the allowable $r$ interval. Take $a,b,c,d$
as the first four cumulative sums of the entries of $\ell$, put
$e=a+d-r_*$, form $J$ by the region rule, and put
$f=2J+q-a-b-c-d-e$. Evaluate $F$ on this tuple and on the tuple
obtained after each legal pair for the type is decreased. The tuple
is used only to evaluate expressions. Its
first five parities are fixed by $\ell$ and $r_*$, and its total
has the prescribed residue relative to $2J$. By \eqref{eq:offset-bounds}
and the normalization above, the local check gives the same result for
this tuple as for every position it represents.

The table gives the enumeration totals. The accompanying
\texttt{check\_six.py} implements exactly these steps using integer
arithmetic and the two unbounded interval endpoints.
\begin{center}
\begin{tabular}{@{}lrrrrrr@{}}
\toprule
Region & $r\le-1$ & $r=0$ & $r=1$ & $r=2$ & $r\ge3$ & Total\\
\midrule
Local checks & 8,802 & 4,377 & 4,401 & 4,401 & 8,793 & 30,774\\
Failures & 0 & 0 & 0 & 0 & 0 & 0\\
\bottomrule
\end{tabular}
\end{center}
Of the $3^6=729$ gap types, three consist of terminal positions.
The other 726 are covered by the table. This proves the recurrence
on the whole domain. Induction on $N$ gives $F=R_6$ and completes
an independent verification of Theorem~\ref{thm:six}.

The supplement also gives an independent encoding of
Lemma~\ref{lem:verification} in quantifier-free linear integer
arithmetic. It uses six unbounded sorted integer variables, the
fifteen possible pairs, and the formula \eqref{eq:six} directly.
The five queries test for a negative value, a nonzero terminal
value, failure to attain $F(x)-1$, a forbidden option at an even
value, and a smaller even option at an odd value. Z3 returned
\texttt{unsat} in all five cases. A separate dynamic program using
only legal moves and Smith's recurrence checked all 906,192 sorted
six-pile positions in which the largest pile has size at most~26.